\documentclass[12pt]{amsart}
\usepackage{amsfonts}
\usepackage{amssymb}
\usepackage{amscd}

\newcommand{\calS}{{\mathcal S}}
\newcommand{\calK}{{\mathcal K}}
\newcommand{\calL}{{\mathcal L}}
\newcommand{\calP}{{\mathcal P}}
\newcommand{\calB}{{\mathcal B}}
\newcommand{\calE}{{\mathcal E}}
\newcommand{\dbH}{{\mathbb H}}
\newcommand{\dbA}{{\mathbb A}}

\newtheorem*{theoremsn}{Theorem}
\newtheorem{theorem}{Theorem}[section]
\newtheorem{lemma}[theorem]{Lemma}

\newtheorem{corollary}[theorem]{Corollary}
\theoremstyle{definition}
\newtheorem{definition}[theorem]{Definition}
\newtheorem{example}[theorem]{Example}

\theoremstyle{remark}
\newtheorem{remark}[theorem]{Remark}
\numberwithin{equation}{section}

\begin{document}

\keywords{Mapping Class Group, Lower Algebraic K-Theory,
Farrell-Jones Isomorphism Conjecture, Strongly Poly-Free Group,
Configuration Space} \subjclass{20F36, 19A31, 19B28, 19D35, 19J10}

\begin{abstract}
We show that the Fibered Isomorphism Conjecture of T. Farrell and
L. Jones holds for various mapping class groups.  In many cases,
we explicitly calculate the lower algebraic $K$-groups, showing
that they do not always vanish.
\end{abstract}

\title{Algebraic K-Theory of Mapping Class Groups}

\author{Ethan Berkove}
\address{Lafayette College\\
 Department of Mathematics \\
 Easton, PA 18042}
\email{berkovee@lafayette.edu}

\author{Daniel Juan-Pineda}
\address{Instituto de Matem\'aticas\\
 UNAM Campus Morelia \\
 Apartado Postal 61-3 (Xangari) \\
 Morelia, Michoac\'an \\
 MEXICO 58089 }
\email{daniel@matmor.unam.mx }

\author{Qin Lu}
\address{Lafayette College\\
 Department of Mathematics \\
 Easton, PA 18042}
\email{luq@lafayette.edu}

\subjclass{20F36, 19A31, 19B28, 19D35, 19J10}

\keywords{Mapping Class Group, Lower Algebraic K-Theory,
Farrell-Jones Isomorphism Conjecture, Strongly Poly-Free Group,
Configuration Space, Fixed Point Data}

\maketitle
\section{Introduction}
Let $\Gamma$ be a torsion-free discrete group. It is a well-known
conjecture that the Whitehead group $Wh(\Gamma)$ of $\Gamma$ must
vanish.  A major tool in the pursuit of this conjecture has been
the \textit{Fibered Isomorphism Conjecture} (\textit{FIC}) of T.
Farrell and L. Jones.  The FIC asserts that the algebraic
$K$-theory groups of $\mathbb{Z}\Gamma$ may be computed from the
corresponding algebraic $K$-theory groups of the virtually cyclic
subgroups of $\Gamma$ (see \cite{FJ} or the Appendix for a precise
formulation).  When FIC holds for a torsion-free group, its
Whitehead group vanishes.

FIC has been verified in many instances: for discrete cocompact
subgroups of virtually connected Lie groups by Farrell and Jones
\cite{FJ}; for pure braid groups by Aravinda, Farrell, and Roushon
\cite{AFR}; for braid groups by Farrell and Roushon \cite{FR}; for
finitely generated Fuchsian groups by Berkove, Juan-Pineda, and
Pearson \cite{BJP}; and for Bianchi groups by Berkove,
Juan-Pineda, Farrell, and Pearson \cite{BFJP}.  In this paper we
investigate when FIC holds for various mapping class groups.

The pure mapping class group $\Gamma_{g,r}^i$ is the group of path
components of orientation preserving self-diffeomorphisms of an
orientable surface of genus $g$, with $i$ punctures and $r$
boundary components.  We require that these diffeomorphisms
pointwise fix the punctures and are the identity on the boundary
components. The full mapping class group is defined similarly, but
includes diffeomorphisms which permute the punctures.

The techniques in this paper build on results in \cite{AFR} and
\cite{FR} concerning strongly poly-free groups.  Although the
mapping class groups are not themselves strongly poly-free, they
admit descriptions where strongly poly-free groups figure
prominently.  We prove in this paper that FIC holds for all pure
mapping class groups of genus $0, 1$ and $2$ and for all full
mapping class groups of genus $0$ and $2$.

\begin{theoremsn}{\rm\textbf{\ref{vanishk}.}}
Let $\Gamma$ be a torsion-free subgroup of any mapping class
group, pure or full, for which FIC holds. Then
$\tilde{K}_i(\mathbb{Z}\Gamma)=0$ for all $i\leq 1$.
\end{theoremsn}

When a mapping class group is torsion-free, Theorem \ref{vanishk}
implies that its Whitehead group is trivial.  However, many
mapping class groups contain considerable torsion.  In these
cases, when FIC holds for a group, often an explicate calculation
is possible.  For example, we have a non-vanishing result for
algebraic $K$-groups.  We note that
\[ \tilde{K}_{i}(\mathbb{Z}\Gamma)=
\begin{cases}
Wh(\Gamma), & \textit{if $i=1;$}\\
\tilde{K}_{0}(\mathbb{Z}\Gamma), & \textit{if $i=0;$}\\
K_{i}(\mathbb{Z}\Gamma),  & \textit{if $i<0.$}
\end{cases}
\]

\begin{theoremsn}{\rm\textbf{\ref{novanishk}.}}
Let $\Gamma$ be a pure mapping class group of genus $g=1$.  Then
for all $i\leq 1$, $\tilde{K}_i(\mathbb{Z}\Gamma) = 0$ with two
exceptions: when $\Gamma = \Gamma_1^0$ or $\Gamma_1^1$, there is
one non-vanishing $K$-group, namely $K_{-1}(\mathbb{Z}\Gamma)
=\mathbb{Z}.$
\end{theoremsn}
\noindent The techniques we develop in the paper also allow us to
study fundamental groups of configuration spaces.

This paper is organized as follows.  In section 2, we give the
basic techniques and results that allow us to apply FIC to mapping
class groups.  In section 3 we prove FIC for various pure mapping
class groups.  We develop these techniques further in section 4 to
study the case of full mapping class groups.  We then use the
implications of FIC to perform explicit $K$-theory calculations in
section 5.  In section 6 we mention how our techniques apply to
other groups, particularly the fundamental groups of configuration
spaces. We mention in the appendix the setup for FIC and some
immediate consequences of its validity.

{\bf Acknowledgements.}  This research was supported in part by
grants CONACyT: 40057-F and DGAPA-UNAM: IN104601.  We would like
to thank Fred Cohen, Tom Farrell, and Jon Pakianathan for helpful
comments and suggestions.  The first author was a visitor at
Lehigh University while this paper was written; he would like to
thank the Department of Mathematics at Lehigh for its generosity
and kind hospitality.

\section{An extension of the Fibered Isomorphism Conjecture}
Before we state our extension of the theorem, we introduce some
background material.  We start with two classes of groups that
figure prominently in our analysis.

\begin{definition}
A group is called {\it virtually cyclic} if it contains a cyclic
group of finite index.
\end{definition}

In particular, a finite group is virtually cyclic.  All infinite
virtually cyclic groups contain an infinite cyclic group of finite
index.  Our arguments also involve a particular type of
torsion-free group whose definition appears in a paper by
Aravinda, Farrell and Roushon.

\begin{definition}
\label{spf} \cite{AFR} A discrete group $\Gamma$ is called {\it
strongly poly-free} if there exists a finite filtration by
subgroups $1 = \Gamma_0 \subseteq \Gamma_1 \subseteq \cdots
\subseteq \Gamma_n = \Gamma$ such that the following conditions
are satisfied:

\begin{enumerate}
\item $\Gamma_i$ is normal in $\Gamma$ for each $i$.
\item $\Gamma_{i+1} / \Gamma_i$ is a finitely generated free group for all $i$.
\item For each $\gamma \in \Gamma$ there is a compact surface $F$
and a diffeomorphism $f: F \to F$ such that the induced
homomorphism $f_{\#}$ on $\pi_1(F)$ is equal to $c_{\gamma}$ in
$Out(\pi_1(F))$, where $c_{\gamma}$ is the action of $\gamma$ on
$\Gamma_{i+1} / \Gamma_i$ by conjugation and $\pi_1(F)$ is
identified with $\Gamma_{i+1} / \Gamma_i$ via a suitable
isomorphism.
\end{enumerate}
\end{definition}
The third condition says that the algebraic action of $\gamma$ on
$\Gamma_{i+1} / \Gamma_i$ can be \textsl{geometrically realized.}

\begin{remark} \label{spfproduct}
It is implicit in the Appendix of \cite{FR} that the finite
product of strongly poly-free groups is also strongly poly-free.
\end{remark}

Our results build on a few known theorems.  The first two are
theorems of Farrell and Jones.

\begin{theorem} \label{extend} \cite[Proposition 2.2]{FJ}
Let $p: \Gamma \to G$ be an epimorphism of groups such that FIC is
true for $G$ and $p^{-1}(S)$ where $S$ ranges over all virtually
cyclic subgroups of $G$.  Then FIC is also true for $\Gamma$.
\end{theorem}

\begin{theorem}\label{inherit}\cite[A.8]{FJ}
If FIC holds for a group $\Gamma$, then FIC also holds for all
subgroups of $\Gamma$.
\end{theorem}

We will also use Farrell and Roushon's Main Theorem, which applies
to extensions by strongly poly-free groups.

\begin{theorem} \label{finiteextend}\cite{FR} Let $\Gamma$ be an extension of a
finite group by a strongly poly-free group (the finite group is
the quotient group). Then $\Gamma$ satisfies FIC.
\end{theorem}

We combine these two results in a theorem that allows us to apply
FIC to the punctured mapping class groups.

\begin{theorem} \label{main}
Say that $\Gamma$ fits into an extension $1 \to K \to \Gamma \to G
\to 1,$ where FIC holds for $G$ and $K$ is a finitely generated
free group.  Furthermore, assume for all $t \in G$ of infinite
order, that the action of the lift $\hat{t}$ on $K$ can be
geometrically realized.  Then $\Gamma$ satisfies FIC.
\end{theorem}

\begin{proof} Let $p: \Gamma \to G$ be the epimorphism in the
short exact sequence.  Given any virtually cyclic subgroup $S
\subseteq G$, we will show that FIC holds for $\hat{S} =
p^{-1}(S)$.

For all virtually cyclic subgroups $S$ there is an extension
$$1 \to K \to \hat{S} \stackrel{p}{\to} S \to 1.$$
We consider two cases: when $S$ is a finite group and when it is
not.  In the case where $|S| < \infty$, Theorem \ref{finiteextend}
directly implies that FIC holds for $\hat{S}$. Next, assume that
$S$ is an arbitrary infinite virtually cyclic group, and let $T$
be a cyclic subgroup of finite index in $S$ with generator $t$.
Note that we may choose $T$ so it is normal in $S$.  Pick a
$\hat{t}$ such that $p(\hat{t}) = t$; the subgroup $<K,\hat{t}> \
\subseteq \hat{S}$ is normal in $\hat{S}$ as $T$ is normal in $S$.

We claim that $<K, \hat{t}>$ is a strongly poly-free group.  The
extension
$$ 1 \to K \to <K, \hat{t}> \to \mathbb{Z} \to 1$$
identifies $<K, \hat{t}>$ as the semidirect product $K \rtimes
<\hat{t}>$, which is filtered by $1 = \Gamma_0 \subseteq K
\subseteq K \rtimes <\hat{t}> = \Gamma_2.$ Thus, the first two
conditions for a strongly poly-free group are satisfied. The third
condition is satisfied by our assumption of the action of
$\hat{t}$ on $K$.

We again write $\hat{S}$ as a finite extension, this time with
strongly poly-free kernel $K \rtimes <\hat{t}>$:
$$1 \to K \rtimes <\hat{t}> \to \hat{S} \to F \to 1.$$
Now, Theorem \ref{finiteextend} implies that FIC holds for
$\hat{S}$. As $S$ is arbitrary, Theorem \ref{extend} implies the
result.
\end{proof}

\begin{remark}
The geometric realization condition in Definition \ref{spf} is
necessary. It is not difficult to find algebraic actions which
have no geometric realization (see, for example, Stallings
\cite{S}).  At this point it is not known if FIC holds in general
for cases where geometric realization is not possible.

\end{remark}

We mention one other way that FIC can be extended to new groups.

\begin{lemma} \label{finitekernel}  Given a short exact
sequence
$$1 \to F \to \Gamma \stackrel{p}{\to} Q \to 1$$
with $F$ finite, if FIC holds for $Q$ then it also holds for
$\Gamma$.
\end{lemma}

\begin{proof}  Let $S$ be any
virtually cyclic subgroup of $Q$, and let $T \subseteq S$ be a
cyclic group of finite index, possibly all of $S$.  There is a
short exact sequence
$$1 \to F \to p^{-1}(T) \stackrel{p}{\to} T \to 1.$$
By definition, $p^{-1}(T)$ is a virtually cyclic group, so it
contains a cyclic subgroup $\hat{T}$ of finite index. Furthermore,
$p^{-1}(T)$ is of finite index in $p^{-1}(S)$, which implies that
$p^{-1}(S)$ is virtually cyclic as it also contains $\hat{T}$.
Therefore, FIC holds for $p^{-1}(S)$.  Theorem \ref{extend}
completes the proof as $S$ is arbitrary.
\end{proof}

\section{Pure Mapping Class Groups}

Denote by $S_{g,r}^i$ an orientable surface of genus $g$ with $i$
punctures and $r$ boundary components.  We define
$Diff^+(S_{g,r}^i)$ to be the group of orientation-preserving
diffeomorphisms of $S_{g,r}^i$ which pointwise fix the punctures
and the boundary. Furthermore, let $Iso(S_{g,r}^i)$ be the
subgroup of $Diff^+ (S_{g,r}^i )$ consisting of all
diffeomorphisms that are isotopic to the identity map.  We define
the {\it pure mapping class group}, $\Gamma_{g,r}^i$, as the
quotient $Diff^+ (S_{g,r}^i)/Iso(S_{g,r}^i).$

In a more general construction, one can define $Diff^+(S_{g,r}^i)$
as the set of diffeomorphisms which fix the set of punctures,
possibly permuting them.  The group constructed in this manner is
called the {\it full mapping class group}, which we denote by
$\hat{\Gamma}_{g,r}^i.$ The full mapping class group is a finite
extension of the pure mapping class group
\begin{equation}\label{ses0}
1\to \Gamma_{g,r}^i\to\hat\Gamma_{g,r}^i\to\Sigma_i\to 1,
\end{equation}
where $\Sigma_i$ is the symmetric group on $i$ letters.
 When $i > 0$ and $r =
0$, we refer to the mapping class group as $\Gamma_{g}^i$, and
when both $i$ and $r$ equal $0$, we write $\Gamma_g $ to denote
the {\it unpunctured mapping class group}.

There are two short exact sequences \cite{H} which tie together
different members of the mapping class family.  These sequences
allow us to extend FIC from group to group. They are
\begin{eqnarray}
& 1 \rightarrow  \pi_1(S^i_{g,r}) \rightarrow \Gamma_{g,r}^{i+1}
\rightarrow \Gamma_{g,r}^i \rightarrow 1 &\quad \text{ if }
2g + r + i > 2,  \label{ses1} \\
& 1 \rightarrow  \mathbb{Z} \longrightarrow \Gamma_{g,r+1}^{i-1}
\longrightarrow \Gamma_{g,r}^{i} \rightarrow 1 &\quad \text{ if }
2g +2r + i > 2.  \label{ses2}
\end{eqnarray}
Note that these sequences relate mapping class groups of the
\textit{same} genus.

\begin{theorem}
\label{induct} Assume $2g + r + i > 2$ and $r+i>0$.  If FIC holds
for $\Gamma_{g,r}^i$, then it also holds for (1)
$\Gamma_{g,r}^{i+1}$; and for (2) $\Gamma_{g,r+1}^{i-1}$ when $i>
0$.
\end{theorem}

\begin{proof} Case (1): Since $2g + r + i>2$, we can use Sequence \ref{ses1}. As
$r+i>0$, $\pi_1(S^i_{g,r})$ is a finitely generated  free group,
so Theorem \ref{main} implies the result once we show the
appropriate geometric realizability condition. This follows as
elements of mapping class groups are intrinsically geometric.
Specifically, take $g\in\Gamma_{g,r}^i$ any element of infinite
order, and let $\hat{g}$ be any lift. As $\hat{g}\in
\Gamma_{g,r}^{i+1}$, it can be represented by a
self-diffeomorphism of $S_{g,r}^{i+1}$. By ignoring the
appropriate puncture, it is also a self-diffeomorphism of
$S_{g,r}^i$ inducing the appropriate action.

Case (2): Use Sequence \ref{ses2}. As in Case (1), we look at the
action of $\hat{g}$ on the strongly poly-free group $\mathbb{Z}$.
There are only two possibilities, the trivial and the involution
actions, both of which are geometrically realizable on a cylinder.
\end{proof}

\begin{theorem} \label{case0}
FIC holds for all pure mapping class groups $\Gamma_{0,r}^i$.
\end{theorem}

\begin{proof} It suffices to show that FIC holds for all mapping class
groups $\Gamma_{0,r}^i$ with $r+i \leq 3$. Then we can use Theorem
\ref{induct} and induction on $i$, then $r$, to show that FIC
holds for the rest of the cases.

It is well-known that $\Gamma_0^i \cong 1$ for $ 0 \leq i \leq 3$,
so FIC holds trivially for these groups.  Also, $\Gamma_{0,r}^0
\cong P_{r-1} \times \mathbb{Z}^{r-1}$, where $P_r$ is the pure
braid group on $r$ strings \cite{Bi,H}. Aravinda, Farrell and
Roushon prove in  \cite{AFR} that FIC holds for $P_r \times
\mathbb{Z}^k$ for any $k$.

Applying Case (2) of Theorem \ref{induct} to $\Gamma_0^3$ implies
that FIC holds for $\Gamma_{0,1}^2$.  The same argument implies
that FIC holds for $\Gamma_{0,2}^1$ as $2g + 2r + i = 4$.  For the
case $\Gamma_{0,1}^1$, we start with the observation that
$\Gamma_{0,2}^0 \cong \mathbb{Z}$.  This follows as
$\Gamma_{0,2}^0 \cong P_1 \times \mathbb{Z} \cong \mathbb{Z}.$
Thus, Sequence \ref{ses2} becomes
$$ 1 \rightarrow \mathbb{Z} \rightarrow \mathbb{Z} \rightarrow \Gamma_{0,1}^1
\rightarrow 1,$$ which implies that FIC holds for $\Gamma_{0,1}^1$
as it is finite.  (In fact, as $r>0$, we shall see that this group
is torsion-free, hence trivial.)

\end{proof}

\begin{theorem} \label{case1}
FIC holds for all pure mapping class groups $\Gamma_{1,r}^i$.
\end{theorem}
\begin{proof} We only need consider the cases $r + i \leq 1$ in
order to apply Theorem \ref{induct}.  First, $\Gamma_1^0 \cong
\Gamma_1^1 \cong SL_2(\mathbb{Z}) \cong \mathbb{Z}/4
\ast_{\mathbb{Z}/2} \mathbb{Z}/6$.  As $SL_2(\mathbb{Z})$ acts on
$\mathbb{H}^2$, hyperbolic 2-space, with finite area quotient,
results in \cite{BFJP} imply that FIC holds for this group.

The final case, $\Gamma_{1,1}^0 $, now follows from Case (2) of
Theorem \ref{induct} with $i=1, r=0$.
\end{proof}

We have to prove that FIC holds for mapping class groups one genus
at a time since the value of the genus stays fixed within
Sequences \ref{ses1} and \ref{ses2}. In general, one wants to show
for a fixed value of $g$ that FIC holds for a mapping class group
whose values of $i$ and $r$ are as small as possible.

\begin{remark} \label{stronger} In \cite{R}, Roushon defines a {\it strongly
poly-surface} group. This is similar to a strongly poly-free
group, except that $\Gamma_{i+1} / \Gamma_i$ is isomorphic to the
fundamental group of a surface.  (There is also a technical
condition that is always satisfied when the surface is closed.) In
the proof of Roushon's Main Theorem, he implies that FIC holds for
any extension
$$
1\to K\to\Gamma\to G\to 1,
$$
where $G$ is finite and $K$ is a strongly poly-surface group.
Using this result we get a version of  Theorem \ref{main} with
``finitely generated free group" replaced by ``fundamental group
of a surface." The proof is identical.
\end{remark}

\begin{theorem}
\label{stronginduct} Assume $g \ge 2$. If FIC holds for a mapping
class group $\Gamma_g$, then FIC holds for $\Gamma_{g,r}^i$, for
all $r$ and $i$.
\end{theorem}

\begin{proof} Using the unpunctured mapping class group, Sequence \ref{ses1}
becomes
$$ 1 \rightarrow \pi_1(S_g) \rightarrow
\Gamma_g^1 \rightarrow \Gamma_g \rightarrow 1 $$ where $S_g$ is a
compact surface of genus $g$.  As $\pi_1(S_g)$ is the fundamental
group of a surface, the version of Theorem \ref{main} mentioned in
Remark \ref{stronger} implies that FIC holds for $\Gamma_g^1$.
Then inductively apply Theorem \ref{induct}.
\end{proof}

\section{Full mapping class groups}
We would like to show that FIC holds for the hyperelliptic mapping
class groups and the unpunctured mapping class group $\Gamma_2$.
One route to this result is to work with the full mapping class
groups of genus $0$, which are closely related to these objects.
Sequence \ref{ses0} shows that the full mapping class groups
contain the pure mapping class groups as subgroups of finite
index. Farrell and Roushon develop techniques in \cite{FR} that
are appropriate to such cases, and we adapt them to mapping class
groups.

\begin{definition} Given groups $K$ and $Q$, with $Q$ finite, the {\it wreath
product} of $K$ and $Q$ is the group $K\wr Q=K^{|Q|} \rtimes Q$
where $K^{|Q|}$ is the product of $|Q|$  copies of $K$ indexed by
elements of $Q$, and $Q$ acts on $K^{|Q|}$ via the regular action
of $Q$ on $Q$.
\end{definition}

\begin{theorem}\cite{DM} \label{embedding}  Take a sequence
$1\to K \to G \to Q \to 1$ with $Q$ a finite group. Then there is
an injective homomorphism $G \to K \wr Q$.
\end{theorem}

Recall that when FIC holds for a group it also holds for any
subgroup (this is Theorem \ref{inherit}).  In light of Theorem
\ref{embedding} and Sequence \ref{ses0}, to show that FIC holds
for $\hat\Gamma_{0}^i$ it suffices to show that FIC holds for
$\Gamma_{0}^i\wr \Sigma_i$. We will prove a stronger statement.

\begin{theorem} \label{wreath} Let $Q$ be a finite group. Then FIC
holds for $\Gamma_0^i\wr Q$ for all $i$.
\end{theorem}

\begin{proof} We proceed by induction.  We note
that the theorem holds trivially for $i = 0, 1, 2, 3$, as
$\Gamma_0^i \cong 1$ for these cases.

Let $Q$ be any finite group.  Since $\pi_1(S_0^i)$ is a normal
subgroup of $\Gamma_0^{i+1}$, $\pi_1(S_0^i)^{|Q|}$ is a normal
subgroup of $\Gamma_0^{i+1} \wr Q.$  Therefore, there is another
short exact sequence that comes from Sequence \ref{ses1}
$$
1 \to \pi_1(S_0^i)^{|Q|} \to \Gamma_0^{i+1} \wr Q {\buildrel p
\over \longrightarrow} \Gamma_0^i \wr Q \to 1.
$$
By the induction assumption, FIC holds for $\Gamma_0^i \wr Q$.
Now let $S$ be any virtually cyclic subgroup of $\Gamma_0^i\wr Q$.
If we can show that FIC holds for all groups $p^{-1}(S)$, then we
have that FIC holds for $\Gamma_0^{i+1} \wr Q$.

There are two cases to consider, when $|S| < \infty$ and when $S$
is infinite virtually cyclic.  In the finite case, there is a
short exact sequence
$$1 \to \pi_1(S_0^i)^{|Q|} \to p^{-1}(S) \to S \to 1.$$
The group $\pi_1(S_0^i)^{|Q|}$ is strongly poly-free by Remark
\ref{spfproduct}, hence FIC holds for $p^{-1}(S)$  by Theorem
\ref{finiteextend} .

The second case is when $S$ is an infinite group.  Let $t$ be a
generator of a normal infinite cyclic subgroup of finite index in
$S$. The element $t$ lifts to an element $\hat{t} \in
\Gamma_0^{i+1} \wr Q.$  By taking a larger power of $t$ if
necessary, we can assume that $\hat{t} \in
(\Gamma_0^{i+1})^{|Q|}$, that is, it acts component-by-component.
As in the proof of Theorem \ref{main}, there is a new short exact
sequence
$$ 1 \to \pi_1(S_0^i)^{|Q|} \rtimes <\hat{t}> \to p^{-1}(S) \to Q' \to 1$$
where $Q'$ is a finite group.   Notice that $\pi_1(S_0^i)^{|Q|}
\rtimes <\hat{t}>$ is a subgroup of $\prod_{|Q|} \pi_1(S_0^i)
\rtimes <\hat{t}>$, where the semidirect product action is
considered by coordinate. (This is Fact 2.4 from \cite{FR}.)
Therefore, $p^{-1}(S)$ is a subgroup of a group $\Gamma'$ which
fits into an extension
$$ 1 \to \prod_{|Q|} \pi_1(S_0^i) \rtimes
<\hat{t}> \to \Gamma' \to Q' \to 1.$$ Each term in the product,
$\pi_1(S_0^i) \rtimes <\hat{t}>$, is a strongly poly-free group by
Theorem \ref{main}, so by Remark \ref{spfproduct}, $\prod_{|Q|}
\pi_1(S_0^i) \rtimes <\hat{t}>$ is also a  strongly poly-free
group.  At this point, Theorem \ref{finiteextend} implies that FIC
holds for $\Gamma'$; the subgroup inheritance property in Theorem
\ref{inherit} then implies that FIC also holds for $p^{-1}(S)$.
This completes the induction step and the proof.
\end{proof}

\begin{corollary} \label{fullmcgroup} Let $G$ be a group that fits into an
extension $1 \to \Gamma_0^i \to G \to Q \to 1$ with $|Q| <
\infty$.  Then FIC holds for $G$. In particular, FIC holds for the
full mapping class group $ \hat\Gamma_0^{i}$.
\end{corollary}

\begin{proof}
The first claim follows from Theorem \ref{embedding} and Theorem
\ref{wreath} as $G\subseteq \Gamma_0^i\wr Q$. The second claim
follows from Sequence \ref{ses0}.
\end{proof}

Given an unpunctured surface of genus $g>0$, there is special
diffeomorphism, the \textit{hyperelliptic involution}, which acts
as a reflection across all of the holes. The normalizer of this
involution in $\Gamma_g$ is the {\it hyperelliptic mapping class
group} $\triangle_g$. There is a relationship between the
hyperelliptic mapping class groups and the full punctured mapping
class groups of genus $0$ given by the sequence
$$ 1 \to \mathbb{Z}/2 \to \triangle_g \to \hat{\Gamma}_0^{2g+2} \to 1.$$
It is well-known that $\triangle_2 = \Gamma_2$ \cite{Bi}.  In
general, however, $\triangle_g$ is neither normal nor of finite
index in $\Gamma_g$.  For more information on the hyperelliptic
mapping class groups, the reader should consult \cite{BH} or
\cite{Gr}.

\begin{corollary} \label{hyperelliptic}
FIC holds for the hyperelliptic mapping class groups.
\end{corollary}

\begin{proof}  Apply Theorem \ref{finitekernel} to the short exact sequence
$$ 1 \to \mathbb{Z}/2 \to \triangle_g \to \hat{\Gamma}_0^{2g+2} \to 1.$$
\end{proof}

\begin{corollary}
FIC holds for all pure mapping class groups of genus $2$.
\end{corollary}

\begin{proof}
FIC holds for $\Gamma_2 = \triangle_2$.  Thus, by Theorem
\ref{stronginduct}, FIC holds for $\Gamma_{2,r}^i.$
\end{proof}

As in the case of the pure mapping class groups, one must proceed
genus-by-genus to show that FIC holds for the full mapping class
groups.  Theorem \ref{wreath} is valid for any genus, as long as
Sequence \ref{ses1} holds and an ad hoc calculation for the base
case of the induction is performed.

\begin{example}
We will show that FIC holds for all full mapping class groups of
genus 2. In order to use Theorem \ref{wreath}, we must first show
that FIC holds for $\Gamma_2 \wr Q$, where $Q$ is any finite
group.

We start with the short exact sequence
$$
1 \to \mathbb{Z}/2 \to \Gamma_2 \to \hat{\Gamma}_0^{6} \to 1.
$$
As in the proof of Theorem \ref{wreath} we can build another short
exact sequence
$$
1 \to (\mathbb{Z}/2)^{|Q|} \to \Gamma_2\wr Q \to
\hat{\Gamma}_0^{6}\wr Q \to 1
$$
with finite kernel.  Hence, by Lemma \ref{finitekernel}, if FIC
holds for $\hat{\Gamma}_0^{6}\wr Q$ it will also hold for $
\Gamma_2\wr Q $.  We direct our attention to
$\hat{\Gamma}_0^{6}\wr Q$.

We have another short exact sequence based on Sequence \ref{ses0}:
$$
1 \to (\Gamma_{0}^6)^{|Q|} \to \hat\Gamma_0^6\wr Q \to \Sigma_6\wr
Q \to 1.
$$
Theorem \ref{embedding} implies that $\hat\Gamma_0^6\wr Q
\subseteq (\Gamma_{0}^6)^{|Q|} \wr (\Sigma_6\wr Q).$
Furthermore, $(\Gamma_{0}^6)^{|Q|} \wr (\Sigma_6\wr Q) \subseteq
\Gamma_{0}^6 \wr ((\Sigma_6\wr Q) \times Q)$ by Fact 2.4 of
\cite{FR}.  Since the group $((\Sigma_6\wr Q) \times Q)$ is
finite, FIC holds for $\Gamma_{0}^6 \wr ((\Sigma_6\wr Q) \times Q)
$ by Theorem \ref{wreath}.

Therefore, FIC holds for $\hat{\Gamma}_0^{6}\wr Q$ by subgroup
inheritance, which in turn implies that FIC holds for $
\Gamma_2\wr Q $.  We have now completed a base case for an
inductive argument.  However, we cannot  apply Theorem
\ref{wreath}, because in Sequence \ref{ses1}
$$ 1 \to  \pi_1(S_2) \to \Gamma_2^1 \to \Gamma_2 \to 1, $$
$\pi_1(S_2)$ is not a free group.  However, $\pi_1(S_2)$ is a
strongly poly-surface group, and the proof of Theorem \ref{wreath}
is valid for this case by an identical argument.  This shows FIC
is true for $\Gamma_2^1 \wr Q$, at which point Theorem
\ref{wreath} applies directly.
\end{example}

\section{Calculations}
In this section we explore consequences of the validity FIC as
they apply to the calculation of lower algebraic $K$-groups. Our
first theorem is a vanishing result.

\begin{theorem}\label{vanishk}
Let $\Gamma$ be a torsion-free subgroup of any mapping class
group, pure or full, for which FIC holds. Then
$\tilde{K}_i(\mathbb{Z}\Gamma)=0$ for all $i\leq 1$.
\end{theorem}

\begin{proof}
This follows from Remark \ref{vanish} in the Appendix.
\end{proof}

\begin{remark}
From results in the previous sections, groups that satisfy the
hypotheses of the above theorem include
\begin{enumerate}
\item $\Gamma_{g,r}^i$ for $g =0,1$ and $r>0$ as these groups are torsion-free,
\item $\Gamma_{0,r}^0$ as these are products of pure braid groups and
torsion-free abelian groups,
\item any torsion-free subgroup of any pure mapping class group of genus $\leq 2$,
and any torsion-free subgroup of any full mapping class group of
genus $0$ or $2$.
\item $\Gamma_0^i.$ Recall that
$\Gamma_0^0 \cong \Gamma_0^1 \cong \Gamma_0^2 \cong \Gamma_0^3
\cong 1$.  Therefore, $\Gamma_0^i$ is torsion-free for $i > 3$ as
both the kernel and quotient of Sequence \ref{ses1} are
torsion-free.

\end{enumerate}
\end{remark}

When we want to calculate the lower algebraic K-theory of a
mapping class group which contains torsion there is more work to
do.  A calculation of the lower algebraic $K$-theory of a given
mapping class group consists of three steps.  First, we show that
FIC holds for the group.  Second, we classify the group's
virtually cyclic subgroups.  Finally, we perform the $K$-theory
calculation, using a classifying space for a suitable family of
subgroups of the mapping class group.  We note, however, that the
unpunctured mapping class groups may possibly contain subgroups
isomorphic to $\mathbb{Z} \times \mathbb{Z}/p \times
\mathbb{Z}/p$. The algebraic K-theory of these infinite virtually
cyclic groups is infinitely generated, so the lower algebraic
K-theory of the unpunctured mapping class group will be a
difficult object to calculate in these cases! When $i \geq 1$,
$\mathbb{Z}/p \times \mathbb{Z}/p$ is not a subgroup of the
mapping class group by results in \cite{L}, so a complete lower
algebraic K-theory calculation seems more reasonable.

Using structure theorems for the punctured pure mapping class
groups, it is possible to determine which virtually cyclic groups
can appear as subgroups.  There is a helpful description of
infinite virtually cyclic groups due to Maskit in \cite{M}.

\begin{theorem} An infinite virtually cyclic group $G$ fits into one
of the two following short exact sequences, with $F$ a finite
group.
$$ 1 \rightarrow F \rightarrow G \rightarrow \mathbb{Z} \rightarrow
1$$ or
$$ 1 \rightarrow F \rightarrow G \rightarrow D_\infty \rightarrow
1.$$ In the former case, $G$ is the semidirect product $F \rtimes
\mathbb{Z}$.  In the latter case, $G$ is the amalgamated product
$G_1 \ast_F G_2$, where $F$ is an index two subgroup in both $G_1$
and $G_2$.
\end{theorem}

This characterization yields an effective way to classify infinite
virtually cyclic subgroups when $i>0$.  One first classifies a
group's finite subgroups.  Any virtual cyclic subgroup which
contains a finite group will sit inside the finite group's
normalizer.  Therefore, normalizers will be our next object of
study.  It is proven in \cite{L} that $\Gamma_g^1$ has
$p$-periodic cohomology, which implies that all finite groups
contained in $\Gamma_g^1$ are cyclic.  Since the kernel of the
Sequence \ref{ses1} is torsion-free, the same result applies for
$\Gamma_g^i$.

\begin{corollary}
\label{cyclicclassify} Let $G$ be a virtually cyclic subgroup of a
punctured pure mapping class group.  The following is a list of
possible $G$:
\begin{enumerate}
\item $G$ is a finite cyclic group.
\item $G$ is the direct product $\mathbb{Z}/p \times \mathbb{Z}$
with $p$ prime.
\item $G$ is the semidirect product $\mathbb{Z}/n \rtimes \mathbb{Z}$.
The action of the infinite cyclic generator on all prime order
cyclic subgroups of $\mathbb{Z}/n$ is trivial.
\item $G$ is the amalgamated product
$\mathbb{Z}/2n \ast_{\mathbb{Z}/n} \mathbb{Z}/2n \cong
\mathbb{Z}/n \times D_{\infty}$, with $n$ odd.
\end{enumerate}
\end{corollary}

\begin{proof} For finite groups, the result in \cite{L}
applies.  Next assume that $G$ has infinite order.  Take a finite
cyclic subgroup $\mathbb{Z}/n$ of a mapping class group
$\Gamma_g^i$ and consider its normalizer.  It is an important fact
that when $n$ is prime $N_{\Gamma_g^i}(\mathbb{Z}/n) =
C_{\Gamma_g^i}(\mathbb{Z}/n)$, i.e., the normalizer and
centralizer agree.  Therefore, any extension of $\mathbb{Z}$ by a
prime order cyclic group must be a direct product, proving (2). On
the other hand, an extension of $\mathbb{Z}$ by a composite order
cyclic group will split, yielding $\mathbb{Z}/n \rtimes
\mathbb{Z}$.  Extensions of $\mathbb{Z}$ by prime cyclics appear
as subgroups of this semidirect product, implying (3).  To prove
(4), take $G$ of the form $G_1 \ast_{\mathbb{Z}/n} G_2$.  The only
finite subgroup of $\Gamma_g^i$ which contains $\mathbb{Z}/n$ as
an index two subgroup is $\mathbb{Z}/2n$.  Also, $n$ must be odd;
otherwise, $\mathbb{Z}/n \times D_{\infty}$ contains a copy of
$\mathbb{Z}/2 \times \mathbb{Z}/2.$  We note that the groups
$\mathbb{Z}$ and $D_{\infty}$ appear as degenerate examples of
cases (3) and (4).
\end{proof}

We saw that when $g=0$ the lower algebraic $K$-theory for pure
mapping class groups vanishes. The pure mapping class groups with
$g = 1$ provide more interesting and less trivial examples.  We
will use two results to aid in our calculations: Nielsen's
Realization Theorem, which states that an element of
$\Gamma_{g}^i$ of finite order can be realized as a diffeomorphism
of $S_g$ of the same order; and the Riemann Hurwitz equation. The
Riemann Hurwitz equation uses the notion of a {\it singular
point}, which is defined as follows:  Given any element of order
$n$ in a mapping class group $\Gamma_g^i$, represent it by a
diffeomorphism $f$ of $S_g$ of order $n$ which fixes the $i$
punctures.  There is a projection $\pi$ sending $S_g$ to $S_h =
S_g/<f>$ which is a $n$-sheet branched covering. Let $\hat P$ be a
point of $S_g$ that is fixed by some power of $f$.   Specifically,
there is a largest value $k$ such that $f^{n/k}$ fixes $\hat P$.
We say that the point $P=\pi(\hat P)$ in the quotient space $S_h$
is a {\it singular point of order $k$} with respect to the
$\mathbb Z/n$ action.  Such a singular point has $n/k$ preimages,
$f(\hat{P}), f^2(\hat{P}), \ldots, f^{n/k}(\hat{P}) = \hat{P}$, in
$S_g$ which correspond to a single orbit under $<f>$.  Note that a
fixed point is the singular point of order $n$.

Assume that after projection, the branched covering $S_g \to S_h$
has $q$ singular points $\{P_1, \ldots ,P_q\}$, and that the order
of $P_i$ is $k_i$. The Riemann Hurwitz equation connects the genus
of the surface, the genus of its quotient space, the order of $f$,
and the singular point information together:
$$2g -2 = n \left( (2h-2) + \sum_{i=1}^q (1-1/k_i) \right) $$

The orders of the torsion elements contained in the mapping class
groups of genus $1$ are known, and are summarized below.

\begin{lemma}\cite{L1}
\label{torsion} If $\Gamma_1^i$ has $p$-torsion, then $p=2,3$.
\begin{enumerate}
\item $\Gamma_1^1$ has 2,3 torsion.
\item $\Gamma_1^2$ has 2,3 torsion.
\item $\Gamma_1^3$ has 2,3 torsion.
\item $\Gamma_1^4$ has 2 torsion.
\item $\Gamma_1^i$ has no $p$-torsion for $i\ge 5$.
\end{enumerate}
\end{lemma}

We use this torsion information to classify the infinite virtually
cyclic subgroups.

\begin{theorem}
For punctured mapping class groups of genus $1$,
\begin{enumerate}
\item $\Gamma_1^0\cong \Gamma_1^1\cong Sl_{2}(\mathbb Z)\cong \mathbb
Z/4*_{\mathbb Z/2}\mathbb Z/6$ and the possible virtually cyclic
subgroups in $\Gamma_1^1$ are $\{1, \mathbb Z/2, \mathbb Z/4,
\mathbb Z/3, \mathbb Z/6,\mathbb Z, \mathbb Z/2\times \mathbb Z,
\mathbb Z/4 \rtimes \mathbb Z, \mathbb Z/4 \times \mathbb Z \}$
\item The possible virtually cyclic subgroups in $\Gamma_1^2$ are $\{1,\mathbb
Z/2, \mathbb Z/4$, $\mathbb Z/3,\mathbb Z, \mathbb Z/2\times
\mathbb Z, \mathbb Z/4 \rtimes \mathbb Z, \mathbb Z/4 \times
\mathbb Z\}$.
\item The possible virtually cyclic subgroups in $\Gamma_1^3$ are $\{1, \mathbb
Z/2$, $\mathbb Z/3$,$\mathbb Z$, $\mathbb Z/2\times \mathbb Z\}$.
\item The possible virtually cyclic subgroups in $\Gamma_1^4$ are $\{1,\mathbb
Z/2, \mathbb Z,\mathbb Z/2\times \mathbb Z\}$
\item The possible virtually cyclic subgroups in $\Gamma_1^i$ for $i\ge 5$ are
$\{1, \mathbb Z \}$.
\end{enumerate}
\end{theorem}

\begin{proof} We will work with the mapping class group
$\Gamma_1^1$ only, as the other cases are similar and easier. By
Lemma \ref{torsion}, we know that $\Gamma_1^1$ has only $2$ and
$3$ torsion, but we still need to determine the exponent of that
torsion.  Recall that $\mathbb{Z}/p \times \mathbb{Z}/p$ does not
appear as a subgroup of any punctured mapping class group.

{\bf The 2-torsion exponent}:   We claim that $\mathbb Z/8$ is not
a subgroup of $\Gamma_1^1$.  By way of contradiction, assume that
there is a copy of $\mathbb Z/8$ in $\Gamma_1^1$.  By Nielsen's
Realization Theorem, we can lift $\mathbb Z/8$ in
$Diffeo^+(S_1^1)$.  View $\mathbb Z/8$ as acting on $S_1^1$ with
quotient space $S_h$.  This action has at least one fixed point,
namely, the puncture in $S_1^1$. Therefore, the Riemann Hurwitz
equation
$$2g-2=p^3(2h-2)+p^3(1-1/p)a_1+ p^3(1-1/p^2)a_2+p^3(1-1/p^3)a_3$$
must have a non-negative solution $(a_1,a_2,a_3)$, with $a_3 \neq
0$. The value of $a_3$ will give the number of fixed points of the
$\mathbb{Z}/8$ action.  Substitute $g=1$ and $p=2$ into the
Riemann Hurwitz equation to get
$$0=8(2h-2)+4a_1+6a_2+7a_3.$$  If $h \ge 1$, $(a_1,a_2,a_3)$ has
no positive solution over the integers. If $h=0$, we get solutions
(4,0,0) and (1,2,0).  In both cases, $a_3 = 0$, which means there
is no singular point of order 8, contradicting our assumption.

Using a similar argument, we identify $\mathbb{Z}/4$ as a
subgroups of $\Gamma_1^1$.  Nielsen's Realization Theorem and the
Riemann Hurwitz equation generate two possible solutions for the
action of $\mathbb{Z}/4$ on $S_1^1$: $(a_1,a_2)=(4,0)$ and
$(a_1,a_2)=(1,2)$, where $a_1,a_2$ are the number of singular
points of order $2$ and $4$ respectively.  As $a_2 \neq 0$, we
only need consider $(a_1,a_2)=(1,2)$.  The $\mathbb Z/4$ action on
$S_1$ has two fixed points of order $4$ and one singular point of
order 2.

The group $\mathbb{Z}/2$ appears as a subgroup of $\Gamma_1^1$ as
it is a subgroup of $\mathbb{Z}/4$.   By a fixed point data
argument that we do not include here, we determine that there is
one conjugacy class of $\mathbb Z/2$ in $\Gamma_1^1$. The details
of the fixed point data argument can be found in \cite{L1}.

{\bf The 3-torsion exponent}:  We claim that there is no $\mathbb
Z/9$ in $\Gamma_1^1$.  We substitute $g=1$ and $p=3$ into the
Riemann Hurwitz equation to get $$2g-2=p^2(2h-2)+p^2(1-1/p)a_1+
p^2(1-1/p^2)a_2.$$   There is one solution to the equation,
$(a_1,a_2)=(3,0)$, which is geometrically impossible because
$a_2=0$ and we must have at least one fixed point of order $9$.
However, there is a copy of $\mathbb Z/3$ in $\Gamma_1^1$.  The
corresponding solution to the Riemann Hurwitz equation generates
an action on $S_1$ has three fixed points corresponding to
$a_1=3$.

{\bf The mixed torsion cases}: $\mathbb Z/12$ is not a subgroup of
$\Gamma_1^1$.  This follows because the Riemann Hurwitz equation,
\begin{eqnarray*}
2g-2 & = & 12(2h-2)+12(1-1/2)a_1+ 12(1-1/3)a_2 \\
     & + & 12(1-1/4)a_3 + 12(1-1/6)a_4+12(1-1/12)a_5
\end{eqnarray*}
has no solution with $a_5 \neq 0$, which is needed in order to
have at least one fixed point of order 12. On the other hand,
there is a copy of $\mathbb Z/6$ in $\Gamma_1^1$.  The $\mathbb
Z/6$ action on $S_1$ has one fixed point, one singular point of
order 2 and one singular point of order 3  corresponding to the
solution $(a_1,a_2,a_3)=(1,1,1).$

We have now accounted for all the possible subgroups of finite
order. By Corollary \ref{cyclicclassify}, the only possible
infinite virtually cyclic groups are $\{\mathbb Z, D_{\infty},
\mathbb Z/2 \times \mathbb Z,\mathbb Z/3 \times \mathbb Z, \mathbb
Z/4 \rtimes \mathbb Z ,\mathbb Z/4 \times \mathbb Z, \mathbb Z/6
\rtimes \mathbb Z, \mathbb Z/6 \times \mathbb Z,\mathbb
Z/4*_{\mathbb Z/2}\mathbb Z/4, \mathbb Z/6*_{\mathbb Z/3}\mathbb
Z/6\}.$   From the calculations of the two exponent, we know that
there is only one conjugacy class of $\mathbb{Z}/2$ subgroups, and
from the group presentation it is central. Therefore, there is no
copy of $D_{\infty}$ in $\Gamma_1^1$, hence no copies of $\mathbb
Z/4*_{\mathbb Z/2}\mathbb Z/4$, and $\mathbb Z/6*_{\mathbb
Z/3}\mathbb Z/6$ either.  From results in \cite{L1}, the
normalizer of $\mathbb Z/3$ in $\Gamma_1^1$ is $\mathbb Z/6$.
Thus, there is no copy of either $\mathbb Z/3\times \mathbb Z$,
$\mathbb Z/6 \times \mathbb Z$, or $\mathbb Z/6 \rtimes \mathbb Z$
in $\Gamma_1^1$.  The claim follows.
\end{proof}

\begin{remark}
We can actually do a better job of classifying the infinite
virtual cyclic subgroups by using the fact that $\Gamma_1 =
\Gamma_1^1 \cong SL_2(\mathbb{Z}) \cong \mathbb{Z}/4
\ast_{\mathbb{Z}/2} \mathbb{Z}/6$.   By a direct calculation, the
virtually cyclic subgroups of $SL_2(\mathbb{Z})$ are $1,
\mathbb{Z}/2, \mathbb{Z}/3, \mathbb{Z}/4, \mathbb{Z}/6,
\mathbb{Z}$, and $\mathbb{Z} \times \mathbb{Z}/2$.  (For this last
group, if $s$ and $t$ generate the copies of $\mathbb{Z}/4$ and
$\mathbb{Z}/6$, then $<s^2, st> \cong \  \mathbb{Z}/2 \times
\mathbb{Z} $ as $s^2 = t^3$.)  All other direct and semidirect
products should be trivial by results in Fine's book.  This
implies that there are no copies of either $\mathbb Z/4 \rtimes
\mathbb Z$ or $\mathbb Z/4 \times \mathbb Z$ in
$SL_2(\mathbb{Z})$.  Using Sequence \ref{ses1}, these groups are
not contained in $\Gamma_1^2$ either.
\end{remark}

\begin{theorem}\label{novanishk}
Let $\Gamma$ be a mapping class group of genus $g=1$.  Then for
all $i\leq 1$, $\tilde{K}_i(\mathbb{Z}\Gamma) = 0$ with two
exceptions: when $\Gamma = \Gamma_1^0$ or $\Gamma_1^1$, there is
one non-vanishing $K$-group, namely $K_{-1}(\mathbb{Z}\Gamma)
=\mathbb{Z}.$
\end{theorem}

\begin{proof}
Of the list of all possible virtually cyclic subgroups, the only
one with a non-vanishing lower algebraic $K$-group is
$\mathbb{Z}/6$ \cite{LS}. Therefore Remark \ref{vanish} implies
that the only mapping class groups with possibly nonvanishing
lower algebraic $K$-groups are $\Gamma_1^0$ and $\Gamma_1^1$.
These two groups are isomorphic to $SL_2(\mathbb{Z})\cong
\mathbb{Z}/6*_{\mathbb{Z}/2}\mathbb{Z}/4.$ As FIC holds for this
group, an appropriate Mayer-Vietoris argument \cite{MP} proves the
claim.

\end{proof}

\section{other examples}

The results in this paper apply to all groups which can be formed
in stages like the mapping class groups.  Prominent members of
this family include the classical braid groups, braid groups of
surfaces, and fiber-type arrangements.

\begin{definition}
Let $M$ be a manifold without boundary.  The {\it configuration
space}  of $n$ ordered points of $M$ is the space $F(M,n) =
\{(x_1, x_2, \ldots x_n) \in M^n | x_i \neq x_j \text{ if } i \neq
j \}.$
\end{definition}

When $M = \mathbb{R}^2$, the fundamental group of $F(M,n)$ is {\it
Artin's pure braid group} on $n$ strands.  When $M$ is a surface,
the fundamental group of $F(M,n)$ is called the {\it pure braid
group} of $M$ of $n$ strands.  Let $Q_m$ denote $m$ fixed distinct
points in $M$.  There is a theorem due to Fadell and Neuwirth.

\begin{theorem} \cite{FN} \label{config}
If M is a (not necessarily compact) manifold without boundary,
then there is a fibration of spaces $F(M \setminus Q_k, j-k) \to
F(M,j) \to F(M,k)$.
\end{theorem}

It is known \cite{CP} that when $M$ is $\mathbb{R}^2$ or a compact
surface of genus $g > 0$ the space $F(M \setminus Q_m, n)$ is a
$K(\pi, 1)$.  In these cases, the fibration in Theorem
\ref{config} leads to a short exact sequence of groups.  We
consider a special case of the fibration of aspherical spaces in
Theorem \ref{config},
$$F(M \setminus Q_{m+k}, 1) \to F(M \setminus Q_m, k+1) \to
F(M \setminus Q_m, k)$$  where $m$ may be equal to zero.  Notice
that $F(M,1) = M$, so FIC holds for the fundamental group of
$F(M,1)$, whether $M$ is punctured or not.  There is an associated
sequence of homotopy groups
$$1 \to \pi_1(M \setminus Q_{m+k}) \to \pi_1(F(M \setminus Q_m, k+1)) \to
\pi_1(F(M \setminus Q_m, k)) \to 1.$$  The group $\pi_1(M
\setminus Q_{m+k})$ is a free group, and arguments in \cite{AFR}
show that this sequence satisfies the geometric realizability
condition contained in the description of a strongly poly-free
group. Therefore, Theorem \ref{main} implies that FIC holds for
$\pi_1(F(M \setminus Q_m, n))$ where $M$ is a surface as above.

This result is not new.  A configuration space is the complement
of what is known as a {\it fiber-type arrangement}, and Cohen
proves that FIC holds for fiber-type arrangement groups in
\cite{C}. (A similar result is also proved in the Appendix of
\cite{FR}.) We can adapt Cohen's techniques to Theorem \ref{main},
reproducing his result.

Theorem \ref{main} occasionally offers a slight advantage over the
techniques in \cite{AFR} and \cite{C}.  By working with a group
one extension at a time, it is easier to deal with groups that
contain torsion.  As examples, we apply the techniques in this
paper to show that FIC holds for the configuration spaces for the
$2$-sphere and the real projective plane.  Neither of these cases
are considered in \cite{AFR} and \cite{C} because of the torsion
they contain.

We first consider $F(S^2,k)$, the braid space of the $2$-sphere.
From Formula 1.7 in \cite{Bi}, the spaces in the fibration of
Theorem \ref{config} yield a short exact sequence
$$ 1 \to \pi_1(S^2 \setminus Q_{n-1}) \to \pi_1(F(S^2,n)) \to
\pi_1(F(S^2,n-1)) \to 1$$ for $n > 3$.  It is known that
$\pi_1(F(S^2,3)) \cong \mathbb{Z}/2$. When $n = 4$,
$\pi_1(F(S^2,4))$ contains a free group (which is strongly
poly-free) as a normal subgroup with $\mathbb{Z}/2$ quotient, so
by Theorem \ref{finiteextend}, FIC holds for this case too.
Theorem \ref{main} and induction imply that FIC holds for
$\pi_1(F(S^2,n))$ with $n >4$.

We case of $F(\mathbb{P},k)$, the braid space on the real
projective plane, is similar.   By results of Van Buskirk in
\cite{V}, there is a short exact sequence for $n \geq 3$,
$$ 1 \to \pi_1(\mathbb{P} \setminus Q_{n-1}) \to
\pi_1(F(\mathbb{P},n)) \to \pi_1(F(\mathbb{P},n-1)) \to 1,$$ where
$\pi_1(\mathbb{P} \setminus Q_{n-1}) \cong F_{n-1}$ is a free
group on $n-1$ generators.  Van Buskirk proves that
$\pi_1(F(\mathbb{P},2)) \cong Q_8$, the quaternion group with
eight elements.  To start the induction, Theorem
\ref{finiteextend} implies that FIC holds for
$\pi_1(F(\mathbb{P},3))$, as this group has a strongly poly-free
normal subgroup with finite quotient. Induction and the short
exact sequence of homotopy groups then implies that FIC holds for
the rest of these braid spaces.

\section{Appendix}
We recall the Fibered Isomorphism Conjecture formulated in
\cite{FJ}. Let $\calS:TOP\to SPECTRA$ be a covariant homotopy
functor. Let $\mathbf{B}$ be the category of continuous surjective
maps: objects in $\mathbf{B}$ are continuous maps $p:E\to B$,
where $E,B$ are objects in $TOP$, and morphisms between
$p_1:E_1\to P_1$ and $p_2:E_2\to P_2$ consist of continuous maps
$f:E_1\to E_2$ and $g_1:B_1\to B_2$ making the following diagram
commute
\[
\begin{CD}
 E_1@>f>>E_2\\
 @Vp_1VV @Vp_2VV\\
 B_1@>g>>B_2.
 \end{CD}
\]
In this setup, Quinn \cite{Q} constructs a functor from
$\mathbf{B}$ to $\Omega-SPECTRA.$ The value of this spectrum at
$p:E\to B$ is denoted by
$$
\dbH(B;\calS(p)),
$$
and has the property that its value at the object $E\to *$ is
$\calS(E)$. The map of spectra $\dbA$ associated to
\[
\begin{CD}
 E@>id>>E\\
 @VpVV @VpVV\\
 B@>g>>*.
 \end{CD}
\]
is known as \textit{Quinn} assembly map.

Given a discrete group $\Gamma$, let $\calE$ be a universal
$\Gamma$-space for the family of virtually cyclic subgroups of
$\Gamma$ \cite[Appendix]{FJ} and denote by $\calB$ the orbit space
$\calE/\Gamma$. Let $X$ be any free and properly discontinuous
$\Gamma$- space, and $p:X\times_{\Gamma}\calE\to B$ be the map
determined by the projection onto $\calB$. The \textit{Fibered
Isomorphism Conjecture} (FIC) for $\calS$ and $\Gamma$ is the
assertion that
$$
\dbA : \dbH(\calB;\calS(p))\to \calS(X/\Gamma)
$$
is a weak equivalence of spectra. This conjecture was made in
\cite[1.7]{FJ} for the functors $\calS=\calP(),\calK(),$ and
$\calL^{-\infty}$, the pseudoisotopy, algebraic $K$-theory and
$\calL^{-\infty}$-theory functors. In this paper we
 mean  FIC as FIC for the functor $\calS=\calP()$.

\begin{remark}\label{vanish}
It is known from \cite[Lemma 1.4.2]{FJ} that if a group $\Gamma$
is torsion-free and FIC holds for $\Gamma$, then
$\tilde{K}_i(\mathbb{Z}\Gamma)=0$ for $i\leq 1$. (Note that
$Wh(\Gamma)=\tilde{K}_1(\mathbb{Z}\Gamma).$) Moreover, the same
conclusion is true when $\tilde{K}_i(\mathbb{Z}G)=0$ for $i\leq 1$
and for all virtually cyclic subgroups $G$ of $\Gamma$.
\end{remark}

\end{document}